\documentclass[10pt]{amsart}
\usepackage{amsopn}
\usepackage{amsmath, amsthm, amssymb, amscd, amsfonts, xcolor}
\usepackage[latin1]{inputenc}

\usepackage{xcolor}
\usepackage{enumitem}

\usepackage[OT2,OT1]{fontenc}
\newcommand\cyr{%
\renewcommand\rmdefault{wncyr}%
\renewcommand\sfdefault{wncyss}%
\renewcommand\encodingdefault{OT2}%
\normalfont
\selectfont}
\DeclareTextFontCommand{\textcyr}{\cyr}

\usepackage[pdftex]{hyperref}

\theoremstyle{plain}

\newtheorem{teo}{Theorem}[section]

\newtheorem{prop}[teo]{Proposition}
\newtheorem{lema}[teo]{Lemma}

\theoremstyle{definition}

\newtheorem{ejemplo}[teo]{Example}

\newtheorem{nota}[teo]{Remark}

\numberwithin{equation}{section}

\newcommand\Ad{\operatorname{Ad}}

\newcommand\rk{\operatorname{rk}}

\newcommand\Exp{\operatorname{Exp}}

\newcommand\so{\mathfrak{so}}

\makeatletter\renewcommand\theenumi{\@roman\c@enumi}\makeatother

\title{On the index of  symmetric spaces}

\author{J\"{u}rgen Berndt}
\address{King's College London, Department of Mathematics, London WC2R 2LS, United Kingdom}
\email{jurgen.berndt@kcl.ac.uk}
\thanks{}

\author{Carlos Olmos}
\address{Facultad de Matem\'atica, Astronom\'ia y F\'isica, Universidad Nacional de C\'ordoba, 
Ciudad Universitaria, 5000 C\'ordoba, Argentina}
\email{olmos@famaf.unc.edu.ar}
\thanks{This research was supported by Famaf-UNC and CIEM-Conicet.}


\begin {document}

\begin{abstract}
Let $M$ be an irreducible Riemannian symmetric space. The index of $M$ is the minimal codimension of a (non-trivial) totally geodesic submanifold of $M$. We prove that the index is bounded from below by the rank of the symmetric space. We also classify the irreducible Riemannian symmetric spaces whose index is less or equal than $3$.
\end{abstract}

\maketitle

\section {Introduction}

A submanifold $\Sigma$ of a Riemannian manifold $M$ is said to be totally geodesic if every geodesic in $\Sigma$ is also a geodesic in $M$. The existence and classification of totally geodesic submanifolds are two fundamental problems in submanifold geometry. In this paper we are considering totally geodesic submanifolds in irreducible Riemannian symmetric spaces.

The totally geodesic submanifolds in Riemannian symmetric spaces of rank one were classified by Wolf (\cite{Wo}) in 1963. It is remarkable that the classification of totally geodesic submanifolds in Riemannian symmetric spaces of higher rank is a very complicated and essentially unsolved problem. \'{E}lie Cartan already noticed an algebraic characterization of totally geodesic submanifolds in terms of Lie triple systems. Although a Lie triple system is an elementary algebraic object, explicit calculations with them can be tremendously complicated. Using the Lie triple system approach, Klein obtained between 2008-10 in a series of papers (\cite{K1}, \cite{K2}, \cite{K3}, \cite{K4}) the classification of totally geodesic submanifolds in irreducible Riemannian symmetric spaces of rank two. No complete classifications are known for totally geodesic submanifolds in irreducible Riemannian symmetric spaces of rank greater than two.

A rather well-known result states that an irreducible Riemannian symmetric space which admits a totally geodesic hypersurface must be a space of constant curvature. As far as the authors know, the first proof of this fact was given by Iwahori (\cite{Iw}) in 1965. In 1980, Onishchik introduced in \cite{On} the index $i(M)$ of a Riemannian symmetric space $M$ as the minimal codimension of a  totally geodesic submanifold of $M$. Onishchik gave an alternative proof for Iwahori's result and also classified the irreducible Riemannian symmetric spaces with index $2$.

In this paper we present a new approach to the index based on different methods. Our first main result states:
\begin{teo}\label{bound}
Let $M$ be an irreducible Riemannian symmetric space. Then 
\[\rk(M) \leq i(M).\]
\end{teo}
\noindent Thus the index is bounded from below by the rank of the symmetric space. We prove this result  by showing that for every totally geodesic submanifold $\Sigma$ in $M$ there exists a maximal flat in $M$ which intersects $\Sigma$ and is transversal to $\Sigma$ at a point of intersection.

Our second main result is the classification of all irreducible Riemannian symmetric spaces $M$ of noncompact type with $i(M) \leq 3$. For $i(M) = 1$ and $i(M) = 2$ we provide alternative proofs of the classifications by Iwahori and Onishchik, respectively. The classification for $i(M) = 3$ is new. We emphasize that, in contrast to $\rk(M) \in \{1,2\}$, the totally geodesic submanifolds for $\rk(M) = 3$ are not classified yet.  The classification result is:
\begin{teo}\label{classification}
 Let $M$ be an irreducible Riemannian symmetric space of noncompact type.
\begin{itemize}[leftmargin=.3in]
\item[\rm(1)] $i(M) = 1$ if and only if $M$ is isometric to 
\begin{itemize}
\item[\rm(i)] the real hyperbolic space ${\mathbb R}H^k = SO^o_{1,k}/SO_k$, $k \geq 2$.
\end{itemize}
\item[\rm(2)] $i(M) = 2$ if and only if $M$ is isometric to one of the following spaces:
\begin{itemize}
\item[\rm(i)] the complex hyperbolic space ${\mathbb C}H^k = SU_{1,k}/S(U_1U_k)$, $k \geq 2$;
\item[\rm(ii)] the Grassmannian $G_2^*({\mathbb R}^{k+2}) = SO^o_{2,k}/SO_2SO_k$, $k \geq 3$;
\item[\rm(iii)] the symmetric space $SL_3({\mathbb R})/SO_3$.
\end{itemize}
\item[\rm(3)] $i(M) = 3$ if and only if $M$ is isometric to one of the following spaces:
\begin{itemize}
\item[\rm(i)] the Grassmannian $G_3^*({\mathbb R}^{k+3}) = SO^o_{3,k}/SO_3SO_k$, $k \geq 3$;
\item[\rm(ii)] the symmetric space $G^2_2/SO_4$; 
\item[\rm(iii)] the symmetric space $SL_3({\mathbb C})/SU_3$.
\end{itemize}
\end{itemize}
\end{teo}
\noindent Duality between Riemannian symmetric spaces of noncompact type and of compact type preserves totally geodesic submanifolds. Also, if $M$ is an irreducible Riemannian symmetric space of compact type and $\hat{M}$ is its Riemannian universal covering space (which is also a Riemannian symmetric space of compact type), then $i(M) = i(\hat{M})$. 
Therefore Theorem \ref{classification} leads, via duality and covering maps, to the classification of irreducible Riemannian symmetric spaces of compact type whose index is less or equal than $3$.

This paper is organized as follows. In Section 2 we present some preliminaries and basic facts. In Section 3 we investigate the set of maximal flats in an irreducible Riemannian symmetric space $M$. The main result states that for every connected totally geodesic submanifold $\Sigma$ of $M$ and every point $p \in \Sigma$  there exists a maximal flat $F$ of 
$M$ with $p \in F$ such that $F$ is transversal to $\Sigma$ at $p$. This implies Theorem \ref{bound}. In Section 4 we investigate the geometry of Lie triple systems for which the orthogonal complement is also a Lie triple system. If one of the two Lie triple systems is not semisimple, we establish a relation to (extrinsically) symmetric submanifolds of Euclidean spaces and symmetric R-spaces (which are also known as symmetric real flag manifolds). Finally, in Section 5, we prove Theorem \ref{classification}.

\section {Preliminaries and basic facts}

For the general theory and the classification of Riemannian symmetric spaces we refer to \cite{H}.
Let $M$ be an $n$-dimensional irreducible Riemannian symmetric space of noncompact type. As usual we write $M = G/K$, where $G= I^o(M)$ is the connected component of the isometry group $I(M)$ of $M$ containing the identity transformation, $p \in M$, and $K = G_p$ is the isotropy group of $G$ at $p$. Let $\mathfrak g = \mathfrak k \oplus \mathfrak p$ be the corresponding Cartan decomposition of the Lie algebra $\mathfrak g$ of $G$, where $\mathfrak k$ is the Lie algebra of $K$. We identify the subspace $\mathfrak p$ of $\mathfrak g$ with the tangent space $T_pM$ of $M$ at $p$ in the usual way.  Let $B$ be the Killing form of ${\mathfrak g}$ and $\theta$ the Cartan involution on ${\mathfrak g}$ corresponding to the Cartan decomposition $\mathfrak g = \mathfrak k \oplus \mathfrak p$. Then $\langle X , Y \rangle = -B(X,\theta Y)$ is an ${\rm Ad}(K)$-invariant positive definite inner product on ${\mathfrak g}$. The Riemannian metric on $M$ is, up to homothety, induced from this inner product. Since totally geodesic submanifolds are preserved under homotheties, we can assume that the Riemannian metric on $M$ is the one induced by the inner product $\langle \cdot , \cdot \rangle$. We denote by $\nabla$ the Levi Civita connection of $M$.

Let $\Sigma$ be a connected totally geodesic submanifold of $M$ with $1 \leq m = \dim\Sigma < n$. Since $M$ is homogeneous we can assume that $p \in \Sigma$. 
Then $\mathfrak p' = T_p\Sigma$ is a Lie triple system in $\mathfrak p$, that is, $\lbrack \lbrack \mathfrak p',\mathfrak p' \rbrack, \mathfrak p' \rbrack \subset \mathfrak p'$.
Since every connected totally geodesic submanifold of a Riemannian symmetric space can be extended to a complete totally geodesic submanifold, we can also assume that $\Sigma$ is complete.  Then the subspace $\mathfrak g' = \lbrack \mathfrak p',\mathfrak p' \rbrack \oplus \mathfrak p' \subset \mathfrak k \oplus \mathfrak p = \mathfrak g$ is a subalgebra of $\mathfrak g$ and $\Sigma = G'/K'$, where $G'$ is the connected closed subgroup of $G$ with Lie algebra $\mathfrak g'$ and $K' = G'_p$ is the isotropy group of $G'$ at $p$. We denote by $\nabla'$ the Levi Civita connection of $\Sigma$. 

Since $\Sigma$ is connected, complete and totally geodesic in $M$, it is also a Riemannian symmetric space of nonpositive sectional curvature and its de Rham decomposition is of the form
$\Sigma = \Sigma_0 \times \Sigma_1 \times \ldots \times \Sigma_d$, where $\Sigma_0$ is a (possibly $0$-dimensional) Euclidean space and $\Sigma_1,\ldots,\Sigma_d$ are (possibly $0$-dimensional) irreducible Riemannian symmetric spaces of noncompact type. In particular, $\Sigma$ is simply connected and therefore $K'$ is connected.
The transvection group $\hat{G}$ of $\Sigma$ is $\hat{G} = \Sigma_0 \times G^1 \times \ldots \times G^d$, where $\Sigma_0$ acts on itself by translations and $G^i = I^o(\Sigma_i)$ (see Theorem 8.3.12 in \cite{W}). We have $\Sigma = \hat{G}/\hat{K}$ with $\hat{K} = \hat{G}_p$. Note that $\hat{K}$ is connected since $\hat{G}$ is connected and $\Sigma$ is simply connected.

For $g \in G'$ we denote by $g_{\vert \Sigma}$ the restriction of $g$ to $\Sigma$. The maps $G' \to \hat{G}, g \mapsto g_{\vert \Sigma}$ and $K' \to \hat{K}, g \mapsto g_{\vert \Sigma}$ are  local group isomorphisms. In particular, we have ${\mathfrak g}' = \hat{\mathfrak g}$ and ${\mathfrak k}' = \hat{\mathfrak k}$. 

Let ${\mathfrak K}(M)$ be the Lie algebra of Killing fields on $M$.
Every $X \in {\mathfrak g}$ determines a Killing field $X^*$ on $M$ by
$X.q = X^*_q = \frac{d}{dt}_{\vert t = 0} (t \mapsto \Exp(tX).q)$
for all $q \in M$, where $\Exp : {\mathfrak g} \to G$ denotes the exponential mapping. The map ${\mathfrak g} \to {\mathfrak K}(M), X \to X^*$ is a Lie algebra isomorphism and its inverse map is given by $\beta : {\mathfrak K}(M) \to {\mathfrak g} = {\mathfrak k} \oplus {\mathfrak p}, Z \mapsto (\nabla Z)_p + Z_p$. Note that $\lbrack X,Y \rbrack^* = -\lbrack X^*,Y^* \rbrack$ for all $X,Y \in {\mathfrak g}$.

For every vector field $Y$ on $M$ we denote by $\bar{Y}$ the vector field on $\Sigma$ which is obtained by first restricting $Y$ to $\Sigma$ and then projecting it orthogonally onto the tangent bundle $T\Sigma$ of $\Sigma$, that is, $\bar{Y}_q = (Y_q)^T$ for all $q \in \Sigma$, where $(\cdot )^T$ denotes the orthogonal projection from $TM$ onto $T\Sigma$. Since $\Sigma$ is totally geodesic in $M$, we have $\bar{Y} \in {\mathfrak K}(\Sigma)$ for all $Y \in {\mathfrak K}(M)$.

Let $X \in {\mathfrak g}$. Then $Y = X^* \in {\mathfrak K}(M)$ and $\bar{Y} \in {\mathfrak K}(\Sigma)$. First assume that $X \in {\mathfrak p}$. Then the Killing field $Y = X^*$ on $M$ is $\nabla$-parallel at $p$ and hence, since $\Sigma$ is totally geodesic in $M$, the Killing field $\bar{Y}$ on $\Sigma$ is $\nabla'$-parallel at $p$. Next, assume that $X \in {\mathfrak k}$ 
and let $\xi,\eta \in T_p\Sigma_0 \subset {\mathfrak p}' \subset {\mathfrak p}$. Then we have
\begin{eqnarray*}
\langle \nabla ' _{\xi^*_p} \bar Y, \eta^*_p\rangle & = & 
\langle \nabla _{\xi^*_p} X^*, \eta^*_p\rangle = 
\langle [\xi^*,X^*]_p, \eta^*_p\rangle = 
\langle [X,\xi]^*_p, \eta^*_p\rangle \\ & = &
\langle [X,\xi], \eta\rangle = 
\langle X, [\eta , \xi]\rangle = 0,
\end{eqnarray*}
as $T_p\Sigma_0$ is an abelian subspace of ${\mathfrak p}$. Recall that the transvection algebra $\hat{\mathfrak g}$ of $\Sigma$ is given by $\hat{\mathfrak g} = \Sigma_0 \oplus {\mathfrak g}^1 \oplus \ldots \oplus {\mathfrak g}^d$. From the above calculations we therefore conclude:

\begin{lema}\label {projeted-bounded}
Let $X \in {\mathfrak g}$, $Y = X^*$ the induced Killing field on $M$ and $\bar{Y}$ the induced Killing field on $\Sigma$. Then $\nabla'\bar{Y}_p + \bar{Y}_p \in \hat{\mathfrak g}$.
\end{lema}

\noindent Roughly speaking, if we identify the Lie algebras of isometries and Killing fields as described above, this result states that the induced Killing fields on $\Sigma$ are in the transvection algebra of $\Sigma$. 

The Cartan decomposition ${\mathfrak g} = {\mathfrak k} \oplus {\mathfrak p}$ is a reductice decomposition and therefore the isotropy representation of $K$ on $T_pM$ can be naturally identified 
with the adjoint representation of $K$ on $\mathfrak p$. 
Note that the symmetric space $M$ is irreducible if and only if $K$ acts irreducibly on $\mathfrak p$.
For $k\in K$, $X \in {\mathfrak k}$ and $v\in \mathfrak p$ we define
 $k.v = \Ad(k).v$ and $X.v  = [X,v]$. The orbit of $K$ resp.\ $K'$ containing $v$ is denoted by $K.v$ resp.\ $K'.v$.

\begin {lema} \label {mainP} 
 Let  $0\neq v\in \mathfrak p ' = T_p\Sigma \subset  T_pM = \mathfrak p$. Then $\dim ( K' .v) < \dim (K.v)$.
\end {lema}

\begin {proof}
Since $K'\subset K$ we obviously have $\dim ( K' .v) \leq \dim (K.v)$.
Assume that equality holds, that is, $\dim ( K' .v) = \dim (K.v)$. Both $K$ and $K'$ are connected since $M$ and $\Sigma$ are simply connected, and thus we must have 
$K.v = K'.v \subset \mathfrak p '$. Then the linear span 
of $K.v$ is a proper $K$-invariant subspace of $\mathfrak p$. 
This contradicts the irreducibility of the $K$-action on $\mathfrak p$.
\end {proof}

The following result is well-known. 

\begin {lema}\label {well-known1}  Let $I$ be a countable set, $M$ and $M_i$ ($i \in I$) be smooth manifolds with $\dim(M_i) < \dim(M)$, and $f_i : M_i \to M$ be  smooth maps. Then $\bigcup _{i\in I}f_i (M_i) \neq M$. 
\end {lema}

\section { The set of maximal flats}

Let $r = \rk(M)$ be the rank of $M$. A maximal flat in $M$ is an $r$-dimensional Euclidean space which is embedded in $M$ as a totally geodesic submanifold. In this section we investigate the structure of the set of maximal flats in $M$. 

Let $\mathcal F $ be the set of maximal flats in $M$ containing $p$ 
and let $\mathcal A$ be the set
of maximal abelian subspaces of $\mathfrak p \simeq T_pM$. 
The map ${\mathcal F} \to {\mathcal A}, F \mapsto T_pF$ is an isomorphism and therefore provides a natural identification of both sets with each other.
For $v \in {\mathfrak p}$ we denote by $\mathfrak C (v)  =     
\{w \in \mathfrak p : [v, w] =0\}$ the centralizer of $v$ in ${\mathfrak p}$. A vector $v \in {\mathfrak p}$ is called regular if $\mathfrak C (v) \in {\mathcal A}$.

Let  $\mathfrak a \in \mathcal A$ and $v\in \mathfrak a$ 
be a regular vector. Then $K.v$ is a principal orbit of the $K$-action on ${\mathfrak p}$ and
the normal space $\nu _v (K.v)$ of $K.v$ in ${\mathfrak p}$ is equal to 
$\mathfrak a$. More generally, 
for every $w\in \mathfrak p$ the centralizer $\mathfrak C (w)$ is equal to the normal space $\nu _w (K.w)$.  
Since $K$ acts transitively on the set $\mathcal A$, there exists
for every $\mathfrak a' \in \mathcal A$ a point
$u\in K.v$ such that $\mathfrak a' = \nu _u (K.v)$. 

Since $K$ acts transitively on ${\mathcal A}$ we have
$\mathcal A = K/K_{\mathfrak a}$, where $K_{\mathfrak a} 
= \{ k \in K: k.\mathfrak a = \mathfrak a\}$ is the isotropy group of 
$\mathfrak  a$ in $K$. The isotropy group $K_{\mathfrak a}$ is a compact subgroup of K.
This equips ${\mathcal A}$ (and hence ${\mathcal F}$)
with the structure of a smooth manifold. 

The isotropy group $K_v = \{k \in K : k.v = v\}$ of $K$ at $v$ is a normal subgroup of 
$K_{\mathfrak a}$ and 
$ K_{\mathfrak a} ^o \subset K_v \subset K_{\mathfrak a}$. 
The finite group $W = K_{\mathfrak a}/K_v$ is the so-called 
Weyl group of $\mathfrak a$. 
The Weyl group $W$ may be regarded either as a finite subgroup (generated by reflections) of the orthogonal group $O(\mathfrak a)$, or as 
a finite group acting on $K.v$ by diffeomorphisms (not in 
general by isometries, since $W$ acts from the right on $K.v$). 
The Weyl group $W$ acts irreducibly on $\mathfrak a$ if and only if 
$M$ is irreducible.  

The orbit $K.v = K/K_v $ is a covering space of $\mathcal A = 
K/K_{\mathfrak a}$ and the fibers are the orbits of $W$ 
on $K.v$.  Observe that
$$\dim(\mathcal A) = \dim(K.v) =
\dim (M) - \rk(M)= n- r. $$

For  $z\in \mathfrak p$ we define
$\mathcal A _z : = \{\mathfrak a \in \mathcal A 
: z\in \mathfrak a\}$.
Observe that $\mathcal A _0 = \mathcal A$. 
The statements in the next  lemma are  well-known or easy to  prove.  

\begin {lema} \label {well-known2} Let $z\in \mathfrak p$. 
\begin{itemize}[leftmargin=.3in]
\item[\rm(1)] The centralizer $\mathfrak C (z)$ is a Lie triple system in ${\mathfrak p}$ and the corresponding 
totally geodesic submanifold $N^z$ of $M$ splits off a line in the direction of $z$. 
\item[\rm(2)]The image under  the slice representation  of 
$(K_z)^o$ in $\nu _z (K.z) = \mathfrak C (z)$  coincides with 
the isotropy group $K^z$ at $p$ of the transvection group $G^z$ 
of the symmetric space $N^z$. 
\item[\rm(3)] Any element of $\mathcal A _z$ is a maximal abelian 
subspace of
$\mathfrak C (z)$. 
Conversely, any maximal abelian subspace of $\mathfrak C (z)$ belongs 
to  $\mathcal A _z$.
In particular, the rank of $N^z$ is equal to the rank of $M$.
\item[\rm(4)] The isotropy group $K_z$, or equivalently $(K_z)^o$,  
 acts transitively on  $\mathcal A _z$. 
Thus $\mathcal A _z$ is in a natural way a 
smooth manifold of dimension $\dim (N^z) - \rk(M) 
= n - \dim (K.z) - r$.
\item[\rm(5)] Let $\mathfrak a \in \mathcal A$ and 
$H_1, \ldots , H_s \subset {\mathfrak a}$ be the reflection hyperplanes of the
symmetries of 
the Weyl group 
$W$. Define $J(u) = \{ i \in \{1,\ldots,s\} :  u \in H_i\}$ for $u\in \mathfrak a$ ($J(u) =\emptyset$ if and only if  $u$ is regular). 
Then $J(u) = J (u')$ if and only if 
$\mathfrak C (u) = \mathfrak C (u')$.
\item[\rm(6)] Let $\mathbb V _{\mathfrak a} ^u = \bigcap _{j\in J(u)}
H_j$ (if $J(u)= \emptyset$ then 
$\mathbb V _{\mathfrak a} ^u  = \mathfrak a$). Then $\mathbb V _{\mathfrak a} ^u$ 
is the tangent space at $p$ of the Euclidean factor of $N^z$. 
\end{itemize}
\end {lema}

\begin {nota} 
 If $\mathfrak a \in \mathcal A _z$, then the Weyl group of $\mathfrak a$, regarding 
$\mathfrak a$ as a maximal abelian subspace of $\mathfrak C (z) = T_pN^z$, 
does not act 
irreducibly, since $N^z$ splits off a line. 
So this group does not coincide with $W$.
\end {nota}

We now come to the main result of this section.

\begin {teo}\label {mainT}
Let $M$ be an irreducible Riemannian symmetric space, 
$\Sigma$ a connected totally geodesic submanifold of $M$, and $p \in \Sigma$.  
Then there exists a maximal flat $F$ of 
$M$ with $p \in F$ such that $F$ is transversal to $\Sigma$ at $p$, that is, $T_pF \cap T_p\Sigma = \{0\}$.
\end {teo}

\begin {proof}
Using covering maps (for compact type) and duality between symmetric spaces of compact type and of noncompact type, we can assume that $M$ is of noncompact type.  As $\Sigma$ is an open part of a connected complete totally geodesic submanifold, we can also assume that $\Sigma$ is complete. We continue using the notations from above. We have to prove that there exists a maximal abelian 
subspace $\mathfrak a $ of $\mathfrak p$ such that $\mathfrak 
a \cap \mathfrak p ' = \{0\}$. 

Let $\mathfrak a ' $ be a maximal abelian subspace of $\mathfrak p'$. For any  ${\mathfrak a} \in {\mathcal A}$
we have $ \mathfrak a \cap \mathfrak p '\neq \{0\}$ if and only if there 
exists $ k' \in  K'$ such that 
$  k ' .\mathfrak a \cap \mathfrak a '\neq \{0\}$. 
So, for any  $\mathfrak a \in \mathcal A$ 
 with $\mathfrak a \cap \mathfrak p '\neq \{0\}$, 
there exist $z\in \mathfrak a'$, $\tilde {\mathfrak a} \in 
\mathcal A _z$ and 
$ k '\in K'$ with  $ k' .\tilde {\mathfrak a} = 
\mathfrak a $. 
Let $\mathfrak a_0 \in \mathcal A$ be a maximal abelian subspace containing 
$\mathfrak a '$.

By defining $u'\sim z'$ if and only if $J(u') = J (z')$ we get an equivalence relation on $\mathfrak a _0\setminus \{0\}$. We restrict this equivalence relation to $\mathfrak a '\setminus \{0\}$.
There are finitely many equivalence classes on $\mathfrak a '\setminus \{0\}$, say 
$[z_1], \ldots , [z_s]$ with $z_1, \ldots , z_s \in \mathfrak a'\setminus\{0\}$. 
From Lemma \ref {well-known2}(5) we know that 
for every $z\in \mathfrak a ' \setminus \{0\}$ there exists  $i \in \{1,\ldots,s\}$ such that 
$\mathfrak C (z) = \mathfrak C (z_i)$. 
Lemma \ref {well-known2}(3) then implies
$\mathcal A _z = \mathcal A _{z_i}$.
Let 
$f_i: K' \times \mathcal A _{z_i} \to \mathcal A$ be the smooth map 
defined by $f (k',\mathfrak a) = k'.\mathfrak a$. 
Then $\mathfrak a \in \mathcal A$ intersects 
$\mathfrak p '$ non-trivially if and only if $\mathfrak a $ belongs to the 
union over $i=1, \ldots , s$ of the images of $f_i$. 
But the dimension of 
$K' \times \mathcal A _{z_i} $ is in general not smaller 
than the dimension of $\mathcal A$ (in order to apply Lemma \ref {well-known1} 
to conclude that there exists an element of $\mathcal A$ that intersects 
$\mathfrak p '$ trivially). 
We will replace each $f_i$ by a finite 
number of smooth functions $g^1 _i, \ldots , g^{d(i)}_i$ 
which are defined on smooth manifolds
$X^1_i, \ldots , X^{d(i)}_i$ with $\dim(X^j_i) < \dim({\mathcal A})$ and such that the image of $f_i$ coincides with 
the union over $j= 1, \ldots, d(i)$ of the images of $g^j_i$. The 
result then follows by applying Lemma \ref {well-known1}. 

Let us consider the principal fibre bundle 
\[0\to K'_{z_i} \to 
K'\to K'/K'_{z_i} = K'.z_i\to 0\]
for $i=1, \ldots , s$. 
We cover the compact manifold $K'.z_i$ 
with finitely many open sets $U^1_{i} , \ldots , U^{d(i)}_{i}$ such that 
there exist global sections $\psi ^j_i : U_i^j \to K'$, $j= 1, 
\ldots , d(i)$, $i=1, \ldots , s$. We define a smooth function
$g_i^j: U^j_i \times \mathcal A _{z_i}\to \mathcal A$ by 
$g_i^j (u, \mathfrak a) = f_i (\psi ^j_i (u),\mathfrak a)$. 
Since $K'_{z_i}\subset K$, it leaves invariant the normal space 
$\nu _{z_i}(K.z_i)= \mathfrak C (z_i)$. Then, by Lemma \ref {well-known2}(4), 
$K'_{z_i}$ leaves invariant 
the set $\mathcal A _{z_i}$ . 
Then, with $\Omega _i^j = \psi _i^j (U_i^j). K'_{z_i}$ we get 
 $f_i( \Omega _i^j\times \mathcal A _{z_i}) = 
g_i^j (U_i^j\times \mathcal A _{z_i})$. Since $\Omega _i^1, \ldots , 
\Omega _i^{d(i)}$ cover $K'$, we obtain that 
$\text {Im} (f_i) = \bigcup _{j=1}^{d(i)} \text {Im} (g_i^j)$, where 
$\text {Im} $ denotes the image of the map. Using  Lemma \ref {well-known2}(4) and 
Lemma \ref {mainP} we get
\begin{eqnarray*}
\dim (X_i^j) & = & \dim (U_i^j) + \dim (\mathcal A _{z_i})  
= \dim ( K' .z_i) + \dim (\mathcal A _{z_i}) \\
& = & \dim (K' .z_i) + n - \dim (K.z_i) -r < n -r 
= \dim (\mathcal A).
\end{eqnarray*}
It follows that every $\mathfrak a \in \mathcal A$ with 
$\mathfrak a \cap \mathfrak p ' \neq \{0\}$  belongs to  
the union over  $j= 1 , \ldots, d(i)$ and $i=1, \ldots, s$ of the images of $g_i^j : X_i^j \to \mathcal A$. As $\dim (X_i^j)< \dim (\mathcal A)$, we conclude from Lemma \ref {well-known1}
that there exists a maximal abelian subspace $\mathfrak a \in \mathcal A$ with 
$\mathfrak a \cap \mathfrak p ' = \{0\}$. 
\end {proof}

Theorem \ref{bound} is a consequence of  Theorem \ref{mainT}.

\section {Complementary Lie triple systems and  symmetric submanifolds}

A submanifold $S$ of a Euclidean space ${\mathbb R}^m$ is called a symmetric submanifold if for each point $q \in S$ the orthogonal reflection $\tau$ of ${\mathbb R}^m$ in the normal space $\nu_qS$ leaves $S$ invariant. Symmetric submanifolds in Euclidean spaces were classified by Ferus (\cite{Fe}). Examples of symmetric submanifolds are standard embeddings of symmetric R-spaces.  An orbit of the isotropy representation of a semisimple Riemannian symmetric space is called an R-space (or real flag manifold), and if this orbit is in addition a symmetric space then it is called a symmetric R-space (or symmetric real flag manifold). It turns out that these symmetric R-spaces are of relevance in our context.

\begin {lema}\label {e-sym} 
Let $M=G/K$ be an irreducible simply 
connected Riemannian symmetric space with 
$\rk(M)\geq 2$, where $G= I^o(M)$ and $K = G_p$ is the isotropy group of $G$ at $p \in M$. 
Let $\mathfrak g = \mathfrak k \oplus \mathfrak p$ be the corresponding
Cartan decomposition. Let $v \in {\mathfrak p}$ and assume that 
the orbit  $K.v \subset \mathfrak p = T_pM$ 
is a symmetric submanifold of the Euclidean space ${\mathfrak p}$. Then the tangent space
$T_v(K.v)$ and the normal space $\nu _p (K.v)$ are Lie triple systems and 
the abelian part of $\nu _p (K.v)$ coincides with $\mathbb R v$.
\end {lema}

\begin {proof}
Assume that the orbit $N=K.v$ is a symmetric submanifold of ${\mathfrak p}$ and let $\tau$ be the orthogonal reflection of ${\mathfrak p}$ in the normal space $\nu_vN$. 
Then we have $\tau (N) = N$, 
$d_v\tau(x) = -x$ for all $x \in T_vN$ and
$d_v\tau(\xi) = \xi$ for all $\xi \in \nu_vN$. 
Let $R$ be the Riemannian curvature tensor of $M$ at $p$. Then $R$ takes values in 
$\mathfrak k$ ($K$ is regarded, via the isotropy representation, 
as a subgroup of $SO (T_pM)$). Let $\bar K$ be the Lie subgroup 
of $SO (T_pM)$ generated by $K$ and $\tau (K) = \tau K \tau ^{-1}$.
Then $\bar K$ is not transitive on the unit sphere of $T_pM$, since 
$\bar K .v = N$ and $\rk(M) \geq 2$. Observe that 
both $R$ and $\tau (R)$ take  values  in $\bar{\mathfrak k}$.
Then, by Simons' holonomy theorem (see \cite {O} or \cite{S}), $\tau (R)$ is a scalar 
multiple of $R$. However, both $R$ and $\tau (R)$ have the same (negative) 
scalar curvature, and therefore $\tau (R) = R$ (and so $\tau (K) = K$). 
So we have $\tau = d_ph$ for some  isometry
$h\in I(M)$. The normal space $\nu _v N$ coincides with the set of
fixed vectors 
of $d_ph$ and the tangent space $T_vN$ coincides with the set 
of fixed vectors of $d_p(\sigma \circ h)$, where $\sigma$ 
is the geodesic symmetry of $M$ at $p$. This shows that both $\nu _v N$ 
and $T_vN$ are Lie triple systems. 

We now  show that the abelian part $\mathfrak a$ of 
$\nu _v N$ is spanned by $v$. 
We denote by $\nabla^\perp$ the normal connection of $N$ and by $A_\xi$ the shape operator of $N$ with respect to $\xi \in {\mathfrak a}$.
Since $N$ is contained in the sphere with radius $\|v\|$ in ${\mathfrak p}$, $A_v$ is minus the identity.
From Lemma \ref {projeted-bounded} and Theorem 4.1.7 in \cite{BCO} we obtain that 
$(K_v)^o$ acts trivially on $\mathfrak a$, and hence ${\mathbb R}v \subset {\mathfrak a}$.
Every $\xi \in \mathfrak a$ induces a unique $\nabla ^\perp$-parallel 
and (locally defined) $G$-invariant  normal vector field $\tilde \xi$ 
of $N$ with $\tilde {\xi}_p = \xi$. 

Assume that $\dim({\mathfrak a}) \geq 2$ and let $\xi \in {\mathfrak a}$ with $\xi \notin {\mathbb R}v$.
We claim that $A_\xi$  cannot be a multiple of the 
identity. Otherwise, adding to $\xi$ some scalar multiple of $v$ we obtain
a nonzero element $\psi \in \mathfrak a$  with $A_\psi =0$. Then $\tilde {\psi} $ is constant
on $N$ and so $N$ is not a full submanifold of ${\mathfrak p}$, which contradicts the 
irreducibility of the isotropy representation of $M = G/K$. Thus $A_\xi$ is not a multiple 
of the identity. 

Let $\lambda _1 , \ldots ,\lambda _g$ be 
the different eigenvalues of $A_\xi$, $g\geq 2$. We may assume that 
$\lambda _1 >  0$ and put $z=v + \frac{1}{\lambda_1}\xi$. Then $K.z$ is a 
singular orbit of the $K$-action on ${\mathfrak p}$ with $T_z(K.z) \subsetneq T_vN$. If we decompose $T_vN$ orthogonally into 
$T_v N = T_z(K.z)\oplus \mathbb V$, then
$\nu _z(K.z) = \nu _v N \oplus \mathbb V$. Note that 
$\tau (z)= z$, $d_z\tau(x) = -x$ for all $x \in T_v N = T_z(K.z) \oplus \mathbb V$ and
$d_z\tau(x) = x$ for all $x \in \nu_v N$.
The first normal space of $K.z$ coincides with the normal 
space, since  $K.z$ is a full submanifold with constant principal curvatures 
(see \cite {BCO}). Let $\alpha$ be the second fundamental form of 
$K.z$. Then, if $x,y\in T_z(K.z)$ are arbitrary,  
$d_z\tau (\alpha (x, y)) =\alpha 
(d_z\tau (x),d_z\tau (y)) = \alpha (-x, -y) 
= \alpha (x,y)$.
This implies $d_z\tau(x) = x$ for all 
 $x \in \nu _z (K.z)$, which  
contradicts the fact that 
$d_z\tau (x) = -x$ for all $x \in {\mathbb V} \subset \nu _z (K.z)$.
We thus conclude that $\mathfrak a = \mathbb R v$.
\end {proof}

Recall that a Lie triple system in ${\mathfrak p}$ corresponds to a connected complete totally geodesic submanifold in $M$ and that this submanifold is a simply connected Riemannian symmetric space of nonpositive curvature. The Lie triple system is said to be semisimple if this symmetric space has no Euclidean factor.

\begin {prop}\label {semisimple} 
Let $M=G/K$ be an irreducible simply 
connected Riemannian symmetric space with 
$\rk(M)\geq 2$, where $G= I^o(M)$ and $K = G_p$ is the isotropy group of $G$ at $p \in M$. 
Let $\mathfrak g = \mathfrak k \oplus \mathfrak p$ be the corresponding
Cartan decomposition.
Assume that ${\mathfrak p} = {\mathfrak p}' \oplus {\mathfrak p}''$ decomposes orthogonally into two Lie triple systems $\mathfrak p' , \mathfrak p'' \subset \mathfrak p $.
Then $\mathfrak p ''$ is not semisimple if and only if there exists $v \in {\mathfrak p}''$ such that
the orbit $K.v$ is a symmetric submanifold of ${\mathfrak p}$ with $T_v(K.v) = \mathfrak p'$ and $\nu_v(K.v) = p ''$. 
Moreover, if  $\mathfrak p ''$ is not semisimple, then its  abelian part 
is one-dimensional.
\end {prop} 

\begin {proof}
The ''if"-part and the final statement follow from   Lemma \ref {e-sym}. Conversely, assume that
$\mathfrak p ''$ is not semisimple. Let ${\mathfrak a}$ be the abelian part of $\mathfrak p ''$ and
$0 \neq v \in {\mathfrak a}$. Note that 
$[v , \mathfrak p ''] = \{0\}$ and thus $\mathfrak p ''\subset 
\mathfrak C (v)$, where $\mathfrak C (v)$ is the centralizer of $v$ in ${\mathfrak p}$. Since both ${\mathfrak p}'$ and ${\mathfrak p}''$ are Lie triple systems, there exists an involutive isometry $\tau \in I(M)$ such that 
$\tau (p)= p$, $d_p\tau(x) = -x$ for all $x \in {\mathfrak p}'$ and  
$d_p\tau(x) = x$ for all $x \in {\mathfrak p}''$. 
We have that 
$d_p\tau  (K.v) = d_p\tau  ( K .d_p\tau ^{-1}(v))  =
(d_p\tau   K d_p\tau ^{-1}).v = K.v$.
Since the normal space 
$\nu _v (K.v)$ coincides with $\mathfrak C (v)\supset \mathfrak p ''$, we get 
$ d_p\tau (x) = -x$ for all $x \in T_v(K.v)$. If $\alpha $ is the second fundamental 
form of $K.v$ and $x,y \in T_v(K.v)$, 
then $d_p\tau (\alpha (x,y)) = 
\alpha (d_p\tau (x), d_p\tau (y)) = 
\alpha (-x,-y)= \alpha (x,y)$. So $\text {d}_p\tau $ is the identity 
when restricted to the first normal space of $K.v$. Since 
$K$ acts irreducibly on ${\mathfrak p}$ and $v\neq 0$, the orbit $K.v$ is a full submanifold of 
$\mathfrak p$. Then, since $K.v$ is a submanifold with constant 
principal curvatures, the first normal space coincides with the normal 
space. This implies $T_v(K.v) = \mathfrak p '$, $d_p\tau$
is an extrinsic symmetry of $K.v$ at $p$, and $\nu _v (K.v) = 
\mathfrak C (v) =\mathfrak p ''$
 (see the last part of the proof 
of Lemma \ref {e-sym}). 
\end {proof}

\begin {nota}\label {rank} 
Let $N^v$ be the connected complete totally geodesic submanifold of $M$ corresponding to the Lie triple system ${\mathfrak p}''$ and assume that ${\mathfrak p}''$ is not semisimple. From the proof of Proposition \ref{semisimple} we know that ${\mathfrak p}'' = {\mathfrak C}(v)$. From Lemma \ref {well-known2}(3) we obtain that $\rk(N^v) = \rk(M)$.
\end {nota}

\begin {nota} \label {most singular}
Assume that the  abelian part of the 
normal space $\nu _v (K.v)= \mathfrak C (v)$ of the isotropy orbit  
$K.v\subset \mathfrak p$ has dimension at least $2$. Then 
there exists $\xi \in \nu _v (K.v)$ such that $K.(v+\xi)$ is a 
parallel focal orbit of $K.v$. In particular, $\nu _v (K.v)$ 
is properly contained in $\nu _{v+\xi} (K.(v+\xi))$. 
This is a well-known fact that can be proved, for instance, 
using arguments as in the proof of  Lemma \ref {e-sym}.
Note that if the isotropy orbit  $K.v$, $v\neq 0$,  
is most singular, then the abelian part of $\nu _v (K.v)$ must have 
dimension $1$ and so must coincide with $\mathbb R v$. 
\end {nota}

\begin {nota}\label {isotropy-irr} 
We recall here, in the notation of this paper, 
 Corollary 2.8 from
 \cite {ORi}: If $K_v$ acts irreducibly on  
the tangent space $T_v(K.v)$ of the isotropy orbit $K.v$, 
then $K.v$ is a symmetric submanifold of ${\mathfrak p}$. 
We will  use this result in Section 5. 
\end {nota}

\begin {nota}\label {rank2}
If $M=G/K$ is an $n$-dimensional irreducible Riemannian symmetric space of dimension $n\geq3$, 
then $2\text {rk}(M) + 1 \leq n$. In fact, any principal isotropy orbit 
$K.v \subset \mathfrak p$ has dimension $m\geq 2$ and is isoparametric. 
The  normal space 
$\nu _v(K.v)$ is a maximal abelian subspace of $\mathfrak p$. 
Let $\eta _1 , \ldots , \eta _g \in 
\nu _v (K.v)$ be the curvature normals of $K.v$ at $v$ (see \cite {BCO}). 
For $\rk(M)=1$ the assertion is trivial. We thus assume that
$\rk(M)\geq 2$ and hence $ 2\leq g \leq m$.  
The curvature normals generate $\nu _v(K.v)$ since 
$K.v$ is a full isoparametric submanifold. Moreover, the equality
$g=m$ holds if and only if the curvature normals are mutually perpendicular. 
In this case $K.v$ splits as a product of submanifolds, which contradicts 
the irreducibility of the isotropy representation of $M = G/K$. 
Thus we have $g\leq m-1$ and therefore $\rk(M) = \dim (\nu _v (K.v)) \leq m-1$, which implies 
$n = m + \text {rk}(M) \geq 2\text {rk}(M) +1$.
\end {nota}

\section {Symmetric spaces of index at most three}

In this section we classify all irreducible Riemannian symmetric spaces of noncompact type with $i(M) \leq 3$.
Let $M$ be an irreducible Riemannian symmetric space of noncompact type. From Theorem \ref{bound} we know that $\rk(M) \leq i(M)$, and therefore $\rk(M) \leq 3$ if $i(M) \leq 3$. For $\rk(M) = 1$ the totally geodesic submanifolds were classified by Wolf in \cite{Wo}. Recall that the Riemannian symmetric spaces of noncompact type and with rank equal to one are the real hyperbolic space ${\mathbb R}H^k = SO^o_{1,k}/SO_k$ ($k \geq 2$), the complex hyperbolic space ${\mathbb C}H^k = SU_{1,k}/S(U_1U_k)$ ($k \geq 2)$, the quaternionic hyperbolic space ${\mathbb H}H^k = Sp_{1,k}/Sp_1Sp_k$ ($k \geq 2$) and the Cayley hyperbolic plane ${\mathbb O}H^2 = F_4^{-20}/Spin_9$. The following lemma follows easily from Wolf's classification:

\begin{lema}\label{rank1}
$i({\mathbb R}H^k) = 1$; $i({\mathbb C}H^k) = 2$; $i({\mathbb H}H^k) = 4$; $i({\mathbb O}H^2) = 8$ ($k \geq 2$).
\end{lema}

We will need the following general result:

\begin {lema}\label {Simons}
Let $M=G/K$ be an irreducible Riemannian symmetric space with 
$\rk(M)\geq 2$, where $G= I^o(M)$ and $K = G_p$ is the isotropy group of $G$ at $p \in M$. 
Let $\mathfrak g = \mathfrak k \oplus \mathfrak p$ be the corresponding
Cartan decomposition. Let $\Sigma$ be a nonflat 
totally geodesic submanifold of $M$ such that  $p\in \Sigma$. 
Let $G'$ be the connected subgroup of $G$ with Lie algebra
$[\mathfrak p ', \mathfrak p '] \oplus \mathfrak p ' $,
where $\mathfrak p ' = T_p\Sigma \subset \mathfrak p = T_pM$, and $K' = G'_p$.
Then the slice representation of $(K')^o$ on $\nu _p \Sigma$ 
is nontrivial. 
 \end {lema}

\begin {proof} 
Assume that $(K')^o$ acts trivially on $\nu _p \Sigma$ and consider the orthogonal decomposition
  $T_pM = T_p\Sigma \oplus \nu _p \Sigma$.
Let $R$ and $R'$ be the Riemannian curvature tensors of $M$ and $\Sigma$ at $p$, 
respectively. We define the 
algebraic  curvature tensor 
 $\bar R = R' \oplus 0$ on $T_pM$, where $0$ is the null algebraic 
curvature tensor 
on $\nu _p \Sigma$. By construction, the restriction of $\bar{R}$ to $T_p\Sigma$ 
has values in the isotropy algebra $\mathfrak k'$  
 (regarding $K'$ as a subgroup of 
$SO(T_p\Sigma)$ via the isotropy representation). By assumption, the slice representation of $(K')^o$ is trivial, and therefore
$\bar R$ takes values in the full isotropy algebra $\mathfrak k$.
Then, by the Simons holonomy theorem (see \cite {O} or \cite{S}), $\bar R$ is a scalar 
multiple of $R$. The scalar must be nonzero as $R'\neq 0$. This is  a contradiction, 
since $\bar R$ is degenerated and $R$ is not. 
\end {proof}

We have the following consequence of Proposition \ref {semisimple}:

\begin {lema}\label {main3}
Let $M=G/K$ be an irreducible simply 
connected Riemannian symmetric space with 
$\rk(M)\geq 2$, where $G= I^o(M)$ and $K = G_p$ is the isotropy group of $G$ at $p \in M$. 
Let $\mathfrak g = \mathfrak k \oplus \mathfrak p$ be the corresponding
Cartan decomposition.
Assume that ${\mathfrak p} = {\mathfrak p}' \oplus {\mathfrak p}''$ decomposes orthogonally into two Lie triple systems $\mathfrak p' , \mathfrak p'' \subset \mathfrak p $. Moreover, assume that ${\mathfrak p}'$ is not abelian and $\dim({\mathfrak p}'') \leq 3$.
Let $G'$ be the connected subgroup of $G$ with Lie algebra 
${\mathfrak g}' = [\mathfrak p',\mathfrak p'] \oplus \mathfrak p' $.  If $G'$ does not act  with cohomogeneity 
one on $M$, then $\rk(M) =2$ and $\dim 
(\mathfrak p '')= 3$. Moreover, there exists 
$v\in T_pM$ such that $K.v$ is a symmetric submanifold of  $\mathfrak p$ with $T_p(K.p) = \mathfrak p ' $ and $\nu _p(K.p) = \mathfrak p '' $. 
 \end {lema}

\begin {proof} Let $\Sigma = G'.p$ be the totally geodesic submanifold given by the Lie triple system $\mathfrak p '$. Assume that $G'$ does not act with cohomogeneity one on $M$ and let $K' = G'_p$. Then the slice representation of $(K')^o$ on $\mathfrak p'' = 
\nu _p\Sigma$ is not transitive on the unit sphere. According to Lemma \ref {Simons} this slice 
representation is nontrivial, which implies
$\dim (\rho ((K')^o))=1$ and $\dim (\mathfrak p'') = 3$, where 
$\rho: K'\to O (\nu _p\Sigma) $ is the slice representation. In particular, 
$\rho ((K')^o)$ fixes a vector $0 \neq v \in {\mathfrak p}'' = \nu_p\Sigma$ .
Note that $v$ is unique up to a scalar multiple.
Let $G''$ be the Lie subgroup of $G$ with Lie algebra 
${\mathfrak g}'' = [\mathfrak p'',\mathfrak p''] \oplus \mathfrak p'' $, $K'' = G''_p$ and $\Sigma^\perp = G''.p = G''/K''$ the totally geodesic submanifold determined by the Lie triple system ${\mathfrak p}''$.
Since $K''$ leaves 
$\Sigma $ invariant,  $\rho ((K')^o)$ is an ideal 
of $K''$. This implies that $\Sigma^\perp$ splits off the line 
corresponding to ${\mathbb R}v$ and thus $\rk(\Sigma^\perp) = 2$. 
Therefore $\mathfrak p ''$ is not semisimple.
The assertion then follows from Proposition \ref {semisimple} and Remark \ref {rank}. 
 \end {proof}

\begin {ejemplo} \label {Veronese} 
We illustrate Lemma \ref{main3} with an example. Consider the isotropy representation of the symmetric space $M= SL_3(\mathbb R)/SO_3$. The principal orbits are $3$-dimensional since $\dim(M) = 5$ and $\rk(M) = 2$. The singular orbits are Veronese embeddings of the real projective plane ${\mathbb R}P^2$ into ${\mathfrak p} \cong {\mathbb R}^5$, which are known to be  symmetric submanifolds of ${\mathbb R}^5$.
By Lemma \ref {e-sym}, both the tangent space ${\mathfrak p}'$ and the normal space ${\mathfrak p}''$ at a point of a focal orbit are Lie triple systems in  ${\mathfrak p}$. The corresponding totally geodesic submanifolds of $M$ are a $2$-dimensional real hyperbolic plane $\Sigma = {\mathbb R}H^2$ and, perpendicular to it, a totally geodesic $\Sigma^\perp = {\mathbb R} \times {\mathbb R}H^2$. The  isometry group $G' = SL_2({\mathbb R})$ of $\Sigma$ acts with cohomogeneity $2$ on $M$.
\end {ejemplo}

\subsection {Symmetric spaces of index 1}
The classification of irreducible Riemannian symmetric spaces $M$ of noncompact type with $i(M) = 1$ follows immediately from Theorem \ref{bound} and Lemma \ref{rank1}, but uses Wolf's classification of totally geodesic submanifolds in Riemannian symmetric spaces of rank one. 
We would like to give here an alternative and conceptual proof.

\begin {teo}\label {i=1} 
Let $M$ be an $n$-dimensional irreducible Riemannian symmetric space of noncompact type, $n \geq 2$. If $i(M) = 1$, then $M$ is isometric to the real hyperbolic space ${\mathbb R}H^n = SO^o_{1,n}/SO_n$.
\end {teo} 

\begin {proof}
Since every homogeneous Riemannian manifold of dimension $2$ has constant curvature we can assume that $n \geq 3$.
Since $i(M) = 1$, there exists a connected complete  totally geodesic hypersurface $\Sigma$ in $M$.
Let $p\in \Sigma$, $G= I^o(M)$ and $K = G_p$.
We identify $K$, via the isotropy representation at $p$, 
with a compact subgroup of $SO(T_pM)$. 

If  $\Sigma$ is flat, then $\rk(M)\geq n-1$. Since $M$ is irreducible
this implies that the principal orbits of $K$ are full 
isoparametric submanifolds of $T_pM$ of dimension $1$. 
This implies  $n=2$, which is a contradiction. 
Hence we may assume  that $\Sigma$ is not flat. 

Let $G'$ be the connected subgroup of $G$ with Lie algebra
${\mathfrak g}' = [\mathfrak p ', \mathfrak p '] \oplus \mathfrak p ' $,
where $\mathfrak p ' = T_p\Sigma \subset \mathfrak p = T_pM$, and $K' = G'_p$. Note that $K'$ is connected since $G'$ is connected and $\Sigma$ is simply connected.
Since $\nu _p\Sigma$ has dimension $1$, the slice representation 
of $K'$ on $\nu _p\Sigma$ is trivial. It follows from Lemma \ref 
{Simons} that $\rk(M) = 1$. Since both $T_p\Sigma$ and $\nu_p\Sigma$ are Lie triple systems and $M$ is simply connected, the geodesic reflection $\tau_\Sigma$ of $M$ in $\Sigma$ is a well-defined global isometry of $M$. The differential $d_p\tau_\Sigma$ is the orthogonal reflection of the Euclidean space ${\mathfrak p} = T_pM$ in the hyperplane $T_p\Sigma$. Since $\rk(M) = 1$, the isotropy group $K$ acts transitively on the unit sphere in $T_pM$ and hence on all hyperplanes in $T_pM$. As the orthogonal reflections in all hyperplanes of $T_pM$ generate the orthogonal group $O(T_pM)$,  it follows that $O(T_pM) \subset\Ad(I(M)_p)$ and thus $M$ has constant curvature.
 \end {proof}

\subsection {Symmetric spaces of index 2}

Let $M$ be an $n$-dimensional irreducible Riemannian symmetric space of noncompact type, $n \geq 2$, and assume that $i(M) = 2$. For $n =2$ we have $M = {\mathbb R}H^2 = SO^o_{1,2}/SO_2$ and for $n = 3$ we have $M = {\mathbb R}H^3 = SO^o_{1,3}/SO_3$, which both have $i(M) = 1$. We can therefore assume that $n \geq 4$. Let $\Sigma$ be an $(n-2)$-dimensional totally geodesic submanifold of $M$. We can assume that $\Sigma$ is connected and complete. 
From  Theorem \ref {bound} we know that $\rk(M)\leq 2$. 
If $\Sigma$ is flat we have $2 = \rk(M) \geq n-2$ and therefore $n = 4$. However, there are no $4$-dimensional irreducible Riemannian symmetric spaces $M$ of noncompact type with $\rk(M) = 2$. We can therefore assume that $\Sigma$ is not flat.
Let $p \in \Sigma$, $G'$ be the connected subgroup of $G$ with Lie algebra
${\mathfrak g}' = [\mathfrak p ', \mathfrak p '] \oplus \mathfrak p ' $,
where $\mathfrak p ' = T_p\Sigma \subset \mathfrak p = T_pM$, and $K' = G'_p$.

If $\rk(M) = 2$, it follows from Lemma \ref {Simons} that
the slice representation of $K' = (K')^o$ on 
$\nu _p \Sigma$ is nontrivial 
and hence transitive on the unit sphere. Therefore $G'$ acts on $M$ with 
cohomogeneity one and has a totally geodesic singular orbit $\Sigma$. 
The first author and Tamaru classified in \cite {BT}  the cohomogeneity 
one actions on irreducible Riemannian symmetric space of noncompact type with a totally geodesic singular orbit. From this classification we obtain that $M$ is isometric to the noncompact Grassmannian $G_2^*({\mathbb R}^{k+2}) = SO^o_{2,k}/SO_2SO_k$ (where $2k = n \geq 6$ and $\Sigma = G_2^*({\mathbb R}^{k+1}$)), or to $SL_3({\mathbb R})/SO_3$ (where $n = 5$ and $\Sigma = {\mathbb R} \times {\mathbb R}H^2 $).

If $\rk(M) = 1$, we can use Lemma \ref{rank1} to conclude that $M$ is isometric to the complex hyperbolic space ${\mathbb C}H^k = SU_{1,k}/S(U_1U_k)$ (where $2k = n$ and $\Sigma = {\mathbb C}H^{k-1}$). Altogether this finishes the proof of Theorem \ref{classification} for $i(M) = 2$.

\begin{nota}
Theorem \ref{classification} for $i(M) = 2$ was proved by Onishchik (\cite{On}) with different, mainly algebraic, methods. From  Theorem \ref {bound} we have $\rk(M)\leq 2$, and therefore one can alternatively apply the classifications of totally geodesic submanifolds by Wolf (\cite{Wo}, for $\rk(M) = 1$) and Klein (\cite{K1}, \cite{K2}, \cite{K3}, \cite{K4}, for $\rk(M) = 2$).
\end{nota}

\subsection {Symmetric spaces of index 3}
Let $M$ be an $n$-dimensional irreducible Riemannian symmetric space of noncompact type with $i(M) = 3$. From Theorem \ref{bound} we know that $\rk(M) \leq 3$ and then Lemma \ref{rank1} implies that $\rk(M) \in \{2,3\}$. The smallest dimension of an irreducible Riemannian symmetric space with $\rk(M) \geq 2$ is $n = 5$. We can therefore assume that $n \geq 5$.

The following lemma solves the case $i(M)= 3$ when $\Sigma$ is not semisimple. 

\begin {lema}\label {i=3ns} Let $M$ be an $n$-dimensional irreducible Riemannian symmetric space of noncompact type and $\Sigma$ be an  $(n-3)$-dimensional, connected, complete, totally geodesic submanifold of $M$. Assume that $\Sigma$ is maximal in $M$ (that is, $\Sigma$ is not contained in a totally geodesic submanifold $\bar{\Sigma}$ of $M$ with dimension $n-3 < \dim{\bar{\Sigma}} < n$). 
If $\Sigma$ is not semisimple, then one of the following statements holds:
\begin{itemize}[leftmargin=.3in]
\item[\rm(i)] $M = SL_4({\mathbb R})/SO_4 = SO^o_{3,3}/SO_3SO_3$ and $\Sigma = {\mathbb R} \times SL_3({\mathbb R})/SO_3$;
\item[\rm(ii)] $M = SO^o_{2,3}/SO_2SO_3$ and $\Sigma = {\mathbb R} \times {\mathbb R}H^2$.
\end{itemize}
\end {lema}

\begin {proof}
Obviously we have $i(M) \leq 3$ and therefore $\rk(M) \leq 3$ by Theorem \ref{bound}.
Let $p \in \Sigma$, $G= I^o(M)$, $K = G_p$, $\mathfrak g = \mathfrak k \oplus \mathfrak p$ be the corresponding Cartan decomposition and ${\mathfrak p}' = T_p\Sigma \subset {\mathfrak p}$. Since $\Sigma$ is not semisimple, the abelian part ${\mathfrak a}'$ of ${\mathfrak p}'$ has dimension $\dim({\mathfrak a}') \geq 1$. Let $0 \neq v \in {\mathfrak a}'$. Then we have ${\mathfrak p}' \subset {\mathfrak C}(v) = \nu_v(K.v)$. Since ${\mathfrak C}(v)$ is a Lie triple system (see Lemma \ref{well-known2}) and $\Sigma$ is maximal in $M$ we conclude ${\mathfrak p}' = {\mathfrak C}(v)$ and hence $\dim(K.v) = 3$. From Remark \ref {most singular} 
we see that $\dim (\mathfrak a') = 1$. 
Since $K$ acts irreducibly on ${\mathfrak p}$, the orbit $K.v \subset {\mathfrak p}$ is full in ${\mathfrak p}$. Then $K$ acts effectively (by isometries) on 
$K.v$.  Therefore we must have $3 \leq \dim (K)\leq 6$.

If $\dim(K) = 6$, the classification of symmetric spaces implies $M = G_2^2/SO_4$ or $M = SL_4({\mathbb R})/SO_4 = SO^o_{3,3}/SO_3SO_3$.
The orbit $K.v$ has constant sectional curvature.
Moreover, we have $\dim (K_v)=3$, 
which implies that $K_v$ acts transitively on the $2$-dimensional unit sphere of 
$T_v(K.v)$. In particular, $K_v$ acts irreducibly on  $T_v(K.v)$.
Then, by Remark \ref {isotropy-irr}, $K.v$ is a symmetric 
submanifold of $\mathfrak p$, or equivalently, $K.v$ is an irreducible symmetric $R$-space. 
From the classification of irreducible symmetric $R$-spaces (see, for instance, \cite {BCO}) we obtain 
that $K.v = {\mathbb R}P^3$ and $M= SL_4({\mathbb R})/SO_4$. The corresponding totally geodesic submanifold is $\Sigma = {\mathbb R} \times SL_3({\mathbb R})/SO_3$.

There are no irreducible Riemannian symmetric spaces of noncompact type $M = G/K$ with $\dim(K) = 5$.  

If $\dim (K) = 4$,  then the classification of  symmetric spaces implies that $M = SO^o_{2,3}/SO_2SO_3 = Sp_2({\mathbb R})/U_2$, which has rank $2$ and dimension $6$. This implies $\dim(\Sigma) = 3$.
Therefore the orbit $K.v$ is not principal,
$\dim(K_v) = 1$,  and 
a principal orbit of the $K$-action has dimension $4$. Let $K.(v+\xi)$ be a principal orbit, 
where $\xi \in \nu _v(K.v) = \mathfrak p '$. Note that such an orbit must be 
isoparametric. The normal space 
$\nu _{v+\xi}(K.(v+\xi))$ is an abelian subspace of $\mathfrak p$. 
Observe that $\nu _{v+\xi}(K.(v+\xi))\subset \mathfrak p '$. Moreover, 
 $\nu _{v+\xi}(K.(v+\xi))$ has codimension one 
in $\mathfrak p '$. This implies that $\Sigma$ is either flat or 
splits as a product of a real hyperbolic plane and a Euclidean space. 
However, $\Sigma$ cannot be flat since $\dim(\Sigma) = 3$ and $\rk(M) = 2$, and therefore $\Sigma = {\mathbb R} \times {\mathbb R}H^2$.

If $\dim (K)=3$, then the classification of  symmetric spaces implies that $M = SL_3({\mathbb R})/SO_3$, which has rank $2$ and dimension $5$.
The orbit  $K.v$ is therefore a principal orbit of the isotropy representation, which shows that
 $v$ is regular and so
 the normal space $\nu _v(K.v) = \mathfrak p'$  is abelian. This means that $\Sigma$ is a maximal flat in $M$.
However, a maximal flat in $SL_3({\mathbb R})/SO_3$ is not a maximal totally geodesic submanifold, as it is always contained in a totally geodesic ${\mathbb R} \times {\mathbb R}H^2$.
\end {proof}

The following lemma deals with the case when  $\Sigma$ is semisimple 
and the slice representation of the isotropy group of $\Sigma$ is not transitive.

\begin {lema} \label {i=3s}
Let $M$ be an $n$-dimensional irreducible Riemannian symmetric space of noncompact type with $\rk(M) \geq 2$ and let $\Sigma$ be an  $(n-3)$-dimensional, connected, complete, totally geodesic submanifold of $M$. Assume that $\Sigma$ is maximal in $M$. 
If $\Sigma$ is semisimple and the slice representation of $\Sigma$ is not transitive on the unit sphere in $\nu_p\Sigma$, $p \in \Sigma$, then the following statements hold:
\begin{itemize}[leftmargin=.3in]
\item[\rm(i)] The fixed point set of the slice representation is one-dimensional.
\item[\rm(ii)] The totally geodesic submanifold $\Sigma$ is an Hermitian symmetric space.
\item[\rm(iii)] If $n \geq 6$, the normal space $\nu _p \Sigma$ is a Lie triple system which is not 
semisimple. 
\end{itemize}
\end {lema}

\begin {proof} 
Let $G= I^o(M)$, $K = G_p$ and
$\mathfrak g = \mathfrak k \oplus \mathfrak p$ be the corresponding
Cartan decomposition.  
Let $G'$ be the connected subgroup of $G$ with Lie algebra 
${\mathfrak g}' = [\mathfrak p ', \mathfrak p '] \oplus \mathfrak p ' $,
where $\mathfrak p ' = T_p\Sigma \subset \mathfrak p = T_pM$, and $K' = G'_p$.
According to Lemma \ref {Simons}, 
 the slice representation $\rho : K' \to  
SO (\nu _p \Sigma)$ cannot 
be trivial. Since $\rho (K')$ is not transitive on the unit sphere in $\nu_p\Sigma$, it is a proper connected subgroup of $SO (\nu _p\Sigma) \simeq SO_3$ and so
$\rho (K') \simeq SO_2$.
This implies (i) since $\dim (\nu _p \Sigma) = 3$.

Let $\Sigma = \Sigma_1 \times \ldots \times \Sigma_d$ be the de Rham decomposition of $\Sigma$. 
Let $G^i$ be the connected closed subgroup of $G$ with Lie algebra 
${\mathfrak g}^i= [\mathfrak p ^i, \mathfrak p ^i] \oplus \mathfrak p ^i$,
where $\mathfrak p ^i = T_p\Sigma_i \subset \mathfrak p = T_pM$, and $K^i= G^i_p$.
The isotropy group $K^i$ acts trivially on $T_p\Sigma_j$ for all $j \neq i$. It thus follows from Lemma \ref{Simons} that $\rho (K^i)$  is nontrivial and hence  $\rho (K^i) = 
\rho (K') \simeq SO_2$. This implies that ${\mathfrak k}^i = {\mathfrak h}^i \oplus \so_2$, where ${\mathfrak h}^i$ is the ideal in ${\mathfrak k}^i$ defined by  ${\mathfrak h}^i = \{X \in {\mathfrak k}^i \mid \lbrack X,Y \rbrack = 0\ {\rm for\ all}\ Y \in \nu_p\Sigma\}$. This shows that 
${\mathfrak k}^i$ is not semisimple and hence $\Sigma _i$ is  Hermitian 
symmetric. This implies (ii).

We now assume that $n \geq 6$. Then $\dim(\Sigma) \geq 3$ and hence $\dim(K') \geq 2$. 
We first consider the case that $\Sigma$ is irreducible, that is, $d = 1$. Then $K' = K^1$ and $\dim({\mathfrak h}^1) > 0$. Let $H^1$ be the connected closed subgroup of $K^1$ with Lie algebra ${\mathfrak h}^1$. Since $H^1$ is a normal subgroup of $K^1$, and since $K^1$ acts irreducibly and almost effectively on $T_p\Sigma$, the adjoint action of $H^1$ on $T_p\Sigma$ cannot fix any nonzero vectors in $T_p\Sigma$. By construction, the adjoint action of $H^1$ fixes each vector in $\nu_p\Sigma$. This implies that $\nu_p\Sigma$ is a Lie triple system.

We next consider the case that $\Sigma $ is not irreducible, that is, $d\geq 2$. Let ${\mathfrak z}^i \simeq \so_2$ be the center of ${\mathfrak k}^i$ and $Z^i \simeq SO_2$ be the connected subgroup of $K^i$ with Lie algebra ${\mathfrak z}^i$. Then, as we have seen above, we have $\rho(Z^i) = \rho(K') \simeq SO_2$ for all $i = 1,\ldots,d$, 
 It is not difficult to show that there exist 
nontrivial $z_i \in Z^i$, $i=1, \ldots , d$, such that 
$\rho(k) \in SO(\nu_p\Sigma)$ is the identity, where $k = z_1\cdot \ldots \cdot z_d \in K'$. 
By construction, $k$  does not fix any nonzero vectors in
 $T_p\Sigma$, and therefore 
$\nu _p \Sigma$ coincides with the set of vectors which are fixed by
$k$, which implies that $\nu _p\Sigma$ is a Lie triple 
system. 

Thus we have proved that ${\mathfrak p}'' = \nu _p \Sigma$ is a Lie triple system of $\mathfrak p$ if $\mathfrak p ' =T_p\Sigma$ is semisimple and $n \geq 6$.
Let $G ''$ be the connected closed subgroup of $G$ with Lie algebra 
${\mathfrak g}'' = [{\mathfrak p}'', {\mathfrak p}''] \oplus {\mathfrak p}'' $
 and $K''= G_p''$. The orbit $\Sigma^\perp = G''.p = G''/K''$ is a totally geodesic submanifold of $M$ with $T_p(\Sigma^\perp) = \nu_p\Sigma$. Note that $K''$ is connected since $\Sigma^\perp$ is simply connected and $G''$ is connected. The isotropy group $K''$ normalizes $K'$ since it leaves ${\mathfrak p}'$ invariant. This implies that 
$\rho (K')$ is a normal subgroup of the restriction of $K''$ to 
$\nu _p \Sigma$. Then the one-dimensional subspace of $\nu_p\Sigma$ spanned by the set of fixed vectors of $\rho (K')$ must be 
invariant under $K''$ and thus be fixed pointwise by $K''$ since $K''$ is connected. This shows that $\nu _ p \Sigma$ is not semisimple. 
\end {proof}

We will now prove Theorem \ref{classification} for $i(M) = 3$.
Let $M$ be an $n$-dimensional irreducible Riemannian symmetric space of noncompact type with $i(M) = 3$ and $\Sigma$ be an $(n-3)$-dimensional, connected, complete, totally geodesic submanifold of $M$. Let $p \in \Sigma$,
$G= I^o(M)$, $K = G_p$ and
$\mathfrak g = \mathfrak k \oplus \mathfrak p$ be the corresponding
Cartan decomposition.  
Let $G'$ be the connected closed subgroup of $G$ with Lie algebra
${\mathfrak g}' = [\mathfrak p ', \mathfrak p '] \oplus \mathfrak p '$,
where $\mathfrak p ' = T_p\Sigma \subset \mathfrak p = T_pM$, and $K' = G'_p$.

From Theorem \ref{bound} and Lemma \ref{rank1} we obtain $\rk(M) \in \{2,3\}$. According to Remark \ref{rank2} we have $2\rk(M) + 1 \leq n$ and therefore $n \geq 5$. There is only one irreducible Riemannian symmetric space of noncompact type whose rank is $2$ and dimension is $5$, namely $SL_3({\mathbb R})/SO_3$, which has index $2$ by part (2) of Theorem \ref{classification}. It follows  that $n \geq 6$.

Recall that $\Sigma$ is a simply connected Riemannian symmetric space of nonpositive curvature.
Let $\Sigma = \Sigma_0 \times \Sigma_1 \times \ldots \times \Sigma_d$ be the de Rham decomposition of $\Sigma$, where $\Sigma_0$ is a (possibly $0$-dimensional) Euclidean space and $\Sigma_1,\ldots,\Sigma_d$ are (possibly $0$-dimensional) irreducible Riemannian symmetric spaces of noncompact type.

We first assume that $\dim(\Sigma_0) > 0$, that is, $\Sigma$ is not semisimple. Then, by Lemma \ref{i=3ns}, $M$ is isometric either to $SL_4({\mathbb R})/SO_4 = SO^o_{3,3}/SO_3SO_3$ and $\Sigma = {\mathbb R} \times SL_3({\mathbb R})/SO_3 $, or to
$SO^o_{2,3}/SO_2SO_3$ and $\Sigma = {\mathbb R} \times {\mathbb R}H^2$. However, by part (2) of Theorem \ref{classification}, the index of $SO^o_{2,3}/SO_2SO_3$ is $2$, and therefore $M = SL_4({\mathbb R})/SO_4 = SO^o_{3,3}/SO_3SO_3$ and $\Sigma = {\mathbb R} \times SL_3({\mathbb R})/SO_3 $.

We now assume that $\dim(\Sigma_0) = 0$, that is, $\Sigma$ is semisimple. We first assume that the slice representation of $K'$ on $\nu_p\Sigma$ is not transitive. From Lemma \ref{i=3s} we see that the $3$-dimensional normal space ${\mathfrak p}'' = \nu_p\Sigma$ is a Lie triple system in ${\mathfrak p}$ which splits of a one-dimensional abelian factor ${\mathfrak a}''$. The corresponding connected complete totally geodesic submanifold $\Sigma''$ of $M$ is therefore isometric to ${\mathbb R} \times {\mathbb R}H^2$. Let $0 \neq v \in {\mathfrak a}''$. From Proposition \ref{semisimple} we obtain that the isotropy orbit $K.v\subset \mathfrak p '$ is a symmetric submanifold of ${\mathfrak p}$ with $T_v(K.v) = \mathfrak p'$ and $\nu_v(K.v) = p ''$. From Lemma \ref{i=3s} we know that $\Sigma$ is Hermitian symmetric. Altogether this implies that $K.v$ is an irreducible symmetric R-space which is a Hermitian symmetric space. Since we are in the noncompact case this means that $K.v$ is either an irreducible symmetric R-space of Hermitian type, or an irreducible symmetric R-space of non-Hermitian type whose universal covering space is Hermitian symmetric. The classification of symmetric R-spaces can for instance be found in the appendix of \cite{BCO}. It follows from this list that there is no irreducible symmetric R-space of Hermitian type whose codimension is equal to $3$. The symmetric R-spaces of non-Hermitian type whose universal covering space is Hermitian symmetric are $SO_3/S(O_1O_2) \subset SL_3({\mathbb R})/SO_3$ (which has codimension $3$), $SO_4/S(O_2O_2) \subset SL_4({\mathbb R})/SO_4$ (which has codimension $5$), and $(S^2 \times S^2)/{\mathbb Z}_2 \subset SO^o_{3,3}/SO_3SO_3$ (which has codimension $5$). However, by Theorem \ref{classification}(2), the index of $SL_3({\mathbb R})/SO_3$ is $2$. Altogether it now follows that the slice representation of $K'$ on $\nu_p\Sigma$ is transitive. This implies that $\Sigma$ is a totally geodesic singular orbit of a cohomogeneity one action on $M$. Such singular orbits were classified by the first author and Tamaru in \cite{BT}. From their classification one can easily find those for which the codimension of the singular orbit is $3$ and $\rk(M) \in \{2,3\}$, namely 
\begin{eqnarray*}
\Sigma = SL_3({\mathbb R})/SO_3 & \subset & M =G_2^2/SO_4, \\
\Sigma = SO^o_{3,k-1}/SO_3SO_{k-1} & \subset & M = SO^o_{3,k}/SO_3SO_k \ (k \geq 3), \\ 
\Sigma = {\mathbb R} \times SL_3({\mathbb R})/SO_3 & \subset & M = SL_4({\mathbb R})/SO_4 = SO^o_{3,3}/SO_3SO_3, \\
\Sigma = SL_3({\mathbb R})/SO_3 & \subset & M = SL_3({\mathbb C})/SU_3.
\end{eqnarray*}
This concludes the proof of  Theorem \ref{classification} for $i(M) = 3$.

\begin{nota}
The normal space of $\Sigma = SL_3({\mathbb R})/SO_3 \subset M =G_2^2/SO_4$ is not a Lie triple system. Therefore the assumption in Lemma \ref{i=3s} that 
the slice representation of $\Sigma$ is not transitive on the unit sphere in $\nu_p\Sigma$, $p \in \Sigma$, is essential.
\end{nota}

\begin{nota}
Onishchik calculated in \cite{On} the index of some symmetric spaces. Those with $i(M) \geq 4$ in his list are:
\begin{eqnarray*}
M = SU_{2,k}/S(U_2U_k),\ k \geq 3 & : & i(M) = 4; \\
M = SL_3({\mathbb H})/Sp_3 & : & i(M) = 6; \\
M = Sp_{2,2}/Sp_2Sp_2 & : & i(M) = 6; \\
M = Sp_{2,k}/Sp_2Sp_k,\ k \geq 3 & : & i(M) = 8; \\
M = Sp_4({\mathbb R})/U_4 & : & i(M) = 8;\\
M = E_6^{-78}/F_4 & : & i(M) = 10; \\
M = E_6^{-78}/Spin_{10}U_1 & : & i(M) = 12.
\end{eqnarray*}
\end{nota}

\begin{nota}
Consider the symmetric space $M = G/K = SL_{k+1}({\mathbb R})/SO_{k+1}$, $k \geq 2$, which has $\rk(M) = k$ and $n = \dim(M) =  k(k+3)/2$. Let $0 \neq v \in {\mathfrak p} \cong {\mathbb R}^n$. The orbit $K.v$ has minimal dimension in ${\mathfrak p}$ if and only if $K.v$ is congruent to the Veronese embedding of the real projective space ${\mathbb R}P^k$ into ${\mathbb R}^n$ (see e.g.\ Lemma 8.1 in \cite{ORi}). This ${\mathbb R}P^k$ is an irreducible symmetric R-space of non-Hermitian type and therefore a symmetric submanifold of ${\mathfrak p}$. It follows from Lemma \ref{e-sym} that both $T_v(K.v)$ and $\nu_v(K.v)$ are Lie triple systems in ${\mathfrak p}$. Let $\Sigma$ be the totally geodesic submanifold in $M$ with $T_v\Sigma = \nu_v(K.v)$. Then we have 
\[
k = \rk(M) \leq i(M) \leq {\rm codim}(\Sigma) = \dim({\mathbb R}P^k) = k,
\]
and therefore $i(M) = k$. 
\end{nota}

Motivated by the inequality $\rk(M) \leq i(M)$ it is natural to ask the following questions: 
\begin{itemize}[leftmargin=.3in]
\item[\rm 1.] What are the irreducible Riemannian symmetric spaces $M$ of noncompact type with $\rk(M) = i(M)$? Known examples are:
\begin{eqnarray*}
i(M) = k & : & M = G^*_k({\mathbb R}^{k+n}) = SO^o_{k+n}/SO_kSO_n,\ 1 \leq k \leq n; \\
i(M) = k & : & M = SL_{k+1}({\mathbb R})/SO_{k+1},\ 2 \leq k.
\end{eqnarray*}
\item[\rm 2.] Are there other geometric or algebraic characterizations of the irreducible Riemannian symmetric spaces $M$ of noncompact type with $\rk(M) = i(M)$?
\end{itemize}

\end {document}